\documentclass[12pt]{amsart}
\newtheorem{thm}{Theorem}[section]
\newtheorem{lem}[thm]{Lemma}

\newtheorem{cor}[thm]{Corollary}

\theoremstyle{definition}
\newtheorem{defn}[thm]{Definition}
\theoremstyle{remark}

\newtheorem{exmp}{Example}[section]

\newcommand\pp{X}

\newcommand\mo{\mathcal O}
\newcommand\lra{\longrightarrow}
\newcommand\ra{\rightarrow}

\newcommand\OOmega{\overline{\Omega}^3_X}
\newcommand\Oomega{\overline{\Omega}}

\newcommand\CC{\mathbb C}

\newcommand\PP{\mathbb P}

\newcommand\nn{\mathfrak N}
\newcommand\II{\mathfrak I}

\newcommand\sing{\operatorname{sing}}

\numberwithin{equation}{section}
\newcommand\ppi{\tilde{\pi}}

\newcommand{\VVsect}{\vspace*{0ex}}
\begin{document}
\title{  Defect via differential forms with logarithmic poles}
\author{S\l awomir Cynk,  S\l awomir Rams}

\thanks{Research partially supported by MNSiW grant no. N N201 388834.}

\begin{abstract}
 We  prove
formulae for the Hodge numbers of big resolutions
 of singular hypersurfaces
satisfying  a Bott-type  vanishing
condition.
 \end{abstract}
\maketitle

\vspace*{-1ex}
\section{Introduction}
\label{sec:intro}

The main purpose of this paper is to prove formulae for
the Hodge numbers of big resolutions of certain hypersurfaces with A-D-E singularities.
Let $X$  be a four-dimensional {\sl normal complete variety},
 and let $Y \subset X$ be a hypersurface with
A-D-E singularities such that $\mbox{sing}(X) \cap Y =
\emptyset$. Recall that singularities of $Y$ can be resolved
by consecutive blow-ups of singular and infinitely near singular points of $Y$. We call such a resolution {\sl big}, and
denote it by
$\ppi: \tilde{Y} \ra Y$. We assume that
$$
 h^2(\Oomega^3_{X}) = h^1(\Oomega^3_{X}(Y)) =   h^1({\mathcal O}_{X}(Y + K_{X})) =  0 \mbox{ and }
h^{2}({\mathcal O}_{X}(-Y))= h^{3}({\mathcal O}_{X}(-Y)) = 0,
$$
where $\OOmega :=  {\mathfrak j}_{*}\Omega^3_{\mbox{\tiny reg}(X)}$ and
 ${\mathfrak j}:\mbox{reg}(X) \rightarrow X$ stands for the inclusion. Then the following equalities hold
(see Thm~\ref{ade-osobliwe})
\bgroup\arraycolsep=2pt
\begin{eqnarray*}
 h^{1,1}({\tilde Y})&=& h^{3}(\OOmega) +
 (\chi(\OOmega(Y)) - h^0(\OOmega(Y))) + (\chi({\mathcal O}_{X}(Y + K_{X}))   - h^0({\mathcal O}_{X}(Y + K_{X}))) + \nonumber \\
 & &  + h^1({\mathcal O}_X) - h^4(\OOmega) - (\chi({\mathcal O}_{X}(2Y + K_{X})) - h^0({\mathcal O}_{X}(2Y + K_{X})))
-  h^1({\mathcal O}_{X}(-Y)) \nonumber \\
 & &  +  \mu_Y +  \delta_Y  \, , \nonumber \\
h^{1,2}({\tilde Y})&=&  h^0(\OOmega) + h^2({\mathcal O}_X) + h^0({\mathcal O}_{X}(2Y + K_{X})) -
  h^1(\OOmega)  - h^0({\mathcal O}_{X}(Y + K_{X}))  + \nonumber \\
 & &  -  h^0(\OOmega(Y))  - \mu_Y + \delta_Y \, . \nonumber
\end{eqnarray*}
\egroup
In the above formulae  $\mu_Y$ is the number of singularities and infinitely near singularities of $Y$, whereas
the defect  $\delta_Y$ (see Def.~\ref{def-ideal}) measures how special the position of singular points of $Y$ with respect to
sections of the sheaf ${\mathcal O}_{X}(2 Y + K_{X})$ is. \\
In particular, the assumptions of Thm~\ref{ade-osobliwe} are fulfilled by ample hypersurfaces with A-D-E singularities
in projective (normal) toric fourfolds (Cor.~\ref{toric}), $1$-ample hypersurfaces in complete toric fourfolds
(Ex.~\ref{1-ample}) and $k$-fold solids (Ex.~\ref{nagladkimidzie}). It should be pointed out that we do not require
$Y$ to be quasi-smooth (compare \cite[Thm~10.6]{bat93}). \\
Our interest in the above formulae is justified by their various applications, e.g.
 \cite[Thm~1]{Cynk2}, that has far more restrictive assumptions than our Thm~\ref{ade-osobliwe} (see \eqref{eq-ample}),
 turned out to be useful in the study of
  factoriality  (see e.g. \cite{Chelt}, \cite{digennaro}, \cite{cheltsov2}).

Recall that
in the paper
\cite{Clemens}
 Clemens proves the equality
\begin{equation} \label{eq-clemens}
h^{1,1}(\tilde{Y}_d) = 1 + \mu + \delta \, ,
\end{equation}
where  $Y_d$  is the double cover of  $\PP_3(\CC)$ branched along a degree-$d$ nodal surface $B_d$,  and
 the integer $\delta$, so-called defect,
is defined as the number of dependent conditions
imposed on homogenous forms of degree  $(3/2 \cdot d - 4)$ on $\PP_3$ by
the vanishing in the nodes of $B_d$.

Later, various  generalizations of the above  formula were found
(for a thorough discussion see \cite{Dimca}, \cite[Chapt.~6]{Dimca2}
and \cite{hulek}).
In \cite{Cynk2} the ambient variety $X$ is assumed to be smooth, whereas $Y \subset X$ is a three-dimensional nodal
hypersurface such that
$h^{2}(\Omega^{1}_{X})$, $h^{3}(\Omega^{1}_{X}(-Y))$ vanish and
\begin{equation} \label{eq-ample}
\mbox{the line bundle } {\mathcal O}_{X}(Y) \mbox{  is ample}.
\end{equation}
 If we
 define the defect $\delta$  as the number of dependent equations
that vanishing in the nodes of $Y$ imposes
on the global sections of the bundle $K_{X}(2Y)$, then $h^{1,1}(\tilde{Y})$ is given by the right-hand side of  the formula
 \eqref{eq-clemens}.
 The assumption \eqref{eq-ample}
turns out to be pretty restrictive, for instance ${\mathcal
O}_{X}(Y)$
ceases to be ample as soon  as we blow-up
a point in $X \setminus Y$. In particular,  \cite[Thm~1]{Cynk2} does not imply  \eqref{eq-clemens}.

 In \cite{rams-hab}
 $Y$ is assumed to be  a hypersurface with A-D-E singularities in a {\sl projective normal Cohen-Macaulay fourfold } $X$
 such that $\mbox{sing}(X) \cap Y = \emptyset$ and
$h^{2}(\Oomega^{1}_{X})$, $h^{3}(\Oomega^{1}_{X}(-Y))$ vanish.
Moreover,
 \eqref{eq-ample} is replaced with the following conditions
\begin{equation} \label{eq-prawieample}
h^{i}({\mathcal O}_{X}(-Y))= h^{j}({\mathcal O}_{X}(-2Y))=  h^{j}(\Oomega^{1}_{X}\otimes \mathcal {\mathcal O}_{X}(-Y))=0
\mbox{ for } j = 1, 2,  \mbox{ and  }  i \leq 3.
\end{equation}
Under such assumptions, one has the equality (see  \cite[Thm~4.1]{rams-hab})
$$
  h^{1,1}(\tilde Y)=h^{1}(\Oomega^{1}_{X})+\mu_Y +\delta_Y + h^{3}({\mathcal O}_{X}(-2Y)) \, ,
$$
where
$\delta_Y$
(see  \cite[Lemma~3.3]{rams-hab})
is given by vanishing of sections of ${\mathcal O}_{\pp}(2Y + K_{\pp})$ and their certain directional derivatives
in the singularities of $Y$. In particular, \cite[Thm~4.1]{rams-hab} can be applied to ample hypersurfaces in projective
{\sl simplicial} toric fourfolds.

Both \cite[Thm~1]{Cynk2} and  \cite[Thm~4.1]{rams-hab} are shown by study of cohomologies of the conormal bundle
of the resolution $\tilde{Y}$ in the appropriate blow-up  of $X$.
One expects that the structure of singularities of $X$ should
play no role (see e.g. \cite{vS}) as far as the Hodge
numbers of $\tilde{Y}$ are concerned, but the above approach
depends among others on the use of the Serre duality, so one has to assume that
the singularities of $X$ are mild enough.
Here we apply the properties of  the  Zariski sheaf of germs of
$3$-forms and the Poincar\'e residue map. In this way we
 need neither
\eqref{eq-ample} nor most of its cohomological consequences \eqref{eq-prawieample}. \\
The paper splits in two parts.
In Sect.~\ref{sect-smooth} we  apply certain technical facts from
 \cite{Cynk2} and \cite{rams-hab} to obtain the formulae under the assumption that $X$ is smooth.
Such a result (Thm~\ref{ade}) has
a very simple proof and can be applied in many interesting cases.
Sect.~\ref{sect-singular} is devoted to the main theorem of the paper (Thm~\ref{ade-osobliwe}).

\vspace{1ex}
\noindent
{\em Notations and conventions:} All varieties are defined over
the
base-field $\CC$.
 By a divisor we mean a Weil divisor, and "$\sim$" stands for the linear equivalence.

\section{Smooth ambient variety $X$} \label{sect-smooth}

Let $Y$ be a hypersurface in a smooth four-dimensional projective variety $X$. We assume that all
singularities of $Y$ are A-D-E points. Let  $\ppi : \tilde{Y}  \ra Y$ be the big resolution of $Y$
obtained  as the composition
\begin{equation} \label{eq-prel6}
\ppi =  \sigma_{1} \circ \ldots \circ \sigma_{n} \, : \,
\tilde{Y} \ra Y = : \tilde{Y}^{0} \, ,
\end{equation}
where
 ${\tilde Y} := {\tilde Y}^{n}$ is smooth and
 $ \sigma_{j}:  {\tilde Y}^{j} \ra {\tilde Y}^{j-1}$, for $j = 1, \ldots, n$, is the blow-up with the center
$\sing({\tilde Y}^{j-1}) \neq \emptyset$. Recall that singularities of ${\tilde Y}^{j}$ are isolated double points for each
$j \leq n-1$.  The number of singularities and infinitely near singularities of $Y$ will be denoted by $\mu_Y$.

Let $\tilde{X}^0 := X$ and let  $\tilde{X}^{j}$ stand for
the fourfold obtained from $\tilde{X}^{j-1}$ by blowing it up along  $\sing({\tilde Y}^{j-1})$, $j = 1, \ldots, n$.
We put $\tilde{X} := \tilde{X}^n$. By abuse of notation,
the composition of the blow-ups in question is denoted by  $\ppi : \tilde{X}  \ra X$.

\begin{defn} \label{def-ideal}
Let
$\sum_l k_l E_l:= K_{\tilde{X}/X}$
and let $\pi^{\ast}Y = \tilde{Y} + \sum_l m_l E_l$, where $E_l$ are (reduced) components of the exceptional locus of
 $\ppi : \tilde{X}  \ra X$.
We put
\begin{equation} \label{eq-prel1ipol}
\II_{Y} := \ppi_{\ast}( {\mathcal O}_{{\tilde X}}((k_l - 2 m_l) E_l))
\end{equation}
and define the {\sl defect} of the hypersurface $Y$ as the integer
\begin{equation} \label{eq-def}
\delta_Y =  h^{0}(K_{X}(2 Y) \otimes \II_{Y}) - (h^{0}(K_{X}(2 Y)) - \mu_Y).
\end{equation}
\end{defn}

\noindent
Remarks: 1. The integer given by \eqref{eq-def} coincides with the one defined in \cite[Def~3.1]{rams-hab}
(see also [ibid., (4.2)]). Indeed,
since
\begin{equation} \label{eq-cof-id}
K_{\tilde{X}} + 2   \tilde{Y}  \sim  \ppi^{\ast}(K_X + 2 Y ) +  \sum (k_l - 2 m_l) E_l \, ,
\end{equation}
the projection formula yields the equality
\begin{equation}  \label{eq-prel7}
h^{0}(K_{\tilde{X}}(2 \tilde{Y})) =  h^{0}(K_{X}(2 Y) \otimes \II_{Y}) \, ,
\end{equation}
which, combined with \cite[Lemma~3.2]{rams-hab}, implies our claim. \\
2. Computations in local coordinates with help of \eqref{eq-prel6} show that
$k_l \leq 2 m_l$ and the ideal
$\II_Y$ is given by vanishing of certain directional derivatives in local coordinates
(see the condition II of \cite[Lemma~3.3]{rams-hab}).

In order to render our exposition self-contained,
we collect here several technical facts that will be used in the proof of Thm~\ref{ade}.

From  the Leray spectral sequence and
\eqref{eq-prel6} (see \cite[equalities (2.15),
(2.16)]{rams-hab}), one obtains
\begin{eqnarray*}
& &  h^{i}({\mathcal O}_{X}(-Y)) = h^{i}({\mathcal O}_{{\tilde X}^1}(-\tilde{Y}^1))
 =  \ldots =   h^{i}({\mathcal O}_{\tilde X}(-\tilde{Y}))  \, \quad\text{for all} \quad i,  \label{eq-prel1} \\
& &  h^{i}({\mathcal O}_{X}(-2Y)) = h^{i}({\mathcal O}_{{\tilde X}^1}(-2\tilde{Y}^1)) =
 \ldots =   h^{i}({\mathcal O}_{\tilde X}(-2\tilde{Y})) \, \quad\text{ for } i \leq 2,  \label{eq-prel2}
\end{eqnarray*}
and the  exact sequences  (see \cite[(2.17)]{rams-hab})
\begin{eqnarray*}
0 \lra H^{3}({\mathcal O}_{{\tilde X}^{j-1}}(-2\tilde{Y}^{j-1}))
&\lra& H^{3}({\mathcal O}_{{\tilde X}^j}(-2\tilde{Y}^j))\lra \CC^{\nu_{j-1}} \lra  \\
&\lra& H^{4}({\mathcal O}_{{\tilde X}^{j-1}}(-2\tilde{Y}^{j-1}))\lra H^{4}({\mathcal O}_{{\tilde X}^j}(-2\tilde{Y}^j))\lra 0 \, ,
\end{eqnarray*}
where  $\nu_j$ stands for  the number of points in $\sing({\tilde Y}^{j})$ and $j = 1, \ldots, n$.
Since $\mu_Y = \sum_{j = 0}^{n-1} \nu_j$, the above exact sequences  imply
\begin{equation*} \label{eq-prel3}
\chi({\mathcal O}_{X}(-2Y)) = \chi({\mathcal O}_{{\tilde X}^1}(-2\tilde{Y}^1)) +  \nu_0 = \ldots
 = \chi({\mathcal O}_{\tilde X}(-2\tilde{Y}) ) + \mu_Y \, .
\end{equation*}
Therefore, the Serre duality yields
\begin{equation} \label{eq-prell}
\chi(K_{X}(2Y))  =   \chi(K_{\tilde X}(2\tilde{Y})) + \mu_Y \mbox{ and }  h^i(K_X(Y)) = h^i(K_{\tilde X}(\tilde{Y}))
\mbox{ for every i}.
\end{equation}

By \cite[Lemma~2.3.b]{rams-hab} we have
\begin{equation*}
 h^{i}(\Omega_{X}^1(-Y)) =
 h^{i}(\Omega_{{\tilde X}^{1}}^1(-\tilde{Y}^{1})) = \ldots  =
 h^{i}(\Omega_{{\tilde X}}^1(-\tilde{Y}))  \, \quad\text{for all} \quad i.
\end{equation*}
Thus the Serre duality yields
\begin{equation}  \label{eq-prel4}
 h^{i}(\Omega_X^3(Y)) =  h^{i}(\Omega^3_{{\tilde X}}(\tilde{Y}))  \, \quad\text{for all} \quad i.
\end{equation}

Moreover, if we assume that $h^{2}(\Omega^{1}_{X}) = 0$, then  \cite[Lemma~2.3.a]{rams-hab} gives
\begin{equation*}
h^{1,1}({{\tilde X}^{1}}) = h^{1,1}(X)+\nu_0 \quad \mbox{ and } \quad
h^{1,i}({{\tilde X}^{1}}) = h^{1,i}(X) \quad \mbox{ for } i \not=1.
\end{equation*}
Consequently,  if we assume $h^{1,2}(X) = 0$, we can proceed by induction and apply \cite[Cor.~III.7.13]{hart} to obtain
\begin{equation}  \label{eq-prel5}
h^{3,3}({\tilde X}) = h^{3,3}(X) + \mu_Y  \quad \mbox{ and }
\quad
h^{3,i}({{\tilde X}}) = h^{3,i}(X) \quad \mbox{ for } i \not=3.
\end{equation}

Recall that the sheaf $\Omega^3_{\tilde X}(\log \tilde{Y})$ of
differential $3$-forms with logarithmic poles along $\tilde{Y}$
is defined as
$$
\Gamma(V,\Omega^3_{\tilde X}(\log \tilde{Y})) := \{ \alpha \in  \Gamma(V,\Omega^3_{\tilde X}(\ast \, \tilde{Y}));
\alpha \mbox{ and } \mbox{d}\alpha \mbox{ have at most simple poles along } \tilde{Y} \},
$$
where $V \subset \tilde{X}$ is open and  $\Omega^3_{\tilde X}(\ast \, \tilde{Y}) := \lim\limits_{\substack{\longrightarrow\\k}}
 \, \Omega^3_{\tilde X}(k \tilde{Y})$.

\noindent
We have the folowing exact sequence (see \cite[p.~444]{pet-st})
\begin{equation} \label{poincare}
0 \lra  \Omega^3_{\tilde X} \lra  \Omega^3_{\tilde X}(\log
\tilde{Y}) \lra  \Omega^2_{\tilde Y} \lra  0   \, ,
\end{equation}
and the following resolution of
 the sheaf
 $ \Omega^3_{\tilde X}(\log \tilde{Y}) $
 (see \cite[p.~445]{pet-st})
\begin{equation} \label{resolution}
0 \lra  \Omega^3_{\tilde X}(\log \tilde{Y}) \lra  \Omega^3_{\tilde X}(\tilde{Y}) \lra
 K_{\tilde X}(2 \tilde{Y}) / K_{\tilde X}(\tilde{Y}) \lra 0 \, .
\end{equation}

Now we are in position to prove (compare \cite[Thm~1]{Cynk2},  \cite[Thm~4.1]{rams-hab})

\begin{thm} \label{ade}
Let $X$ be a smooth projective fourfold, and let $Y \subset X$ be a hypersurface with A-D-E singularities. If
$$
h^{2}(\Omega^{1}_{X})= h^{3}(\Omega^{1}_{X}(-Y))=   h^{3}({\mathcal O}_{X}(-Y))=
h^{2}({\mathcal O}_{X}(-Y))= 0 \,,
$$
then
\bgroup\arraycolsep=2pt
\begin{eqnarray*}
 h^{1,1}({\tilde Y})&=& h^{1,1}(X) +  (\chi(\Omega^1_{X}({-Y})) - h^4(\Omega^1_{X}({-Y})))
 - (\chi({\mathcal O}_{X}(-2Y)) - h^4({\mathcal O}_{X}(-2Y))) + \nonumber \\
& &  - 2  h^1({\mathcal O}_{X}(-Y))  +  \mu_Y +  \delta_Y \, ,\nonumber \\
h^{1,2}({\tilde Y})&=&   h^{4,1}(X)  +  h^{0,2}(X) +  h^0(K_X(2Y))   -
  h^{3,1}(X) -  h^4(\Omega^1_{X}({-Y})) -  h^0(K_X(Y)) + \\
& &  - \mu_Y  +   \delta_Y,
\end{eqnarray*}
\egroup
where $\delta_Y$ (resp. $\mu_Y$) is the defect
(resp. the number of singularities and infinitely near singularities) of $Y$.
\end{thm}
\begin{proof}  To simplify our notation we put  $\nn := K_{\tilde X}(2 \tilde{Y}) / K_{\tilde X}(\tilde{Y})$.
From the exact sequence
\begin{equation} \label{eq-iloraz}
0 \lra K_{\tilde{X}}(\tilde{Y}) \lra K_{\tilde{X}}(2 \tilde{Y})  \lra   \nn  \lra 0 \, ,
\end{equation}
we obtain
\begin{eqnarray} \label{eq-chi}
& &  \chi(\nn) =  \chi(K_{\tilde{X}}(2 \tilde{Y})) -  \chi(K_{\tilde{X}}(\tilde{Y}))
\stackrel{\scriptstyle{\eqref{eq-prell}}}{=}  \chi(K_{X}(2Y)) - \mu_Y -
\chi(K_{X}(Y)) \, .
\end{eqnarray}
By \eqref{eq-prell} we have
$h^1(K_{\tilde{X}}(\tilde{Y})) = h^1(K_{X}(Y)) = h^3({\mathcal
O}_{X}(-Y)) = 0$,
so the cohomology sequence associated to \eqref{eq-iloraz}
breaks into shorter exact sequences, which implies
\begin{eqnarray} \label{eq-hoN}
h^0(\nn)  &=&    h^0(K_{\tilde{X}}(2 \tilde{Y})) -   h^0(K_{\tilde{X}}(\tilde{Y}))
 \stackrel{\scriptstyle{\eqref{eq-prel7}}}{=}   h^{0}(K_{X}(2 Y) \otimes
\II_{Y}) - h^0(K_{X}(Y)) =  \nonumber \\
  &\stackrel{\scriptstyle{\eqref{eq-def}}}{=}&  \delta_Y + h^0(K_X(2Y)) - \mu_Y
-  h^0(K_X(Y)). \nonumber
\end{eqnarray}

Observe that  \eqref{eq-prel4} gives $h^{1}(\Omega^{3}_{\tilde{X}}(\tilde{Y}))= 0$.
Thus  the cohomology sequence associated to \eqref{resolution} breaks into shorter exact sequences and yields
\begin{eqnarray}
& &  \sum_{j=0}^{1} (-1)^j \cdot h^j(\Omega^3_{\tilde X}({\log \tilde Y})) =
 h^0(\Omega^3_{\tilde X}({\tilde Y})) -  h^0(\nn)
\stackrel{\scriptstyle{\eqref{eq-prel4}}}{=}
 h^0(\Omega^3_{X}({Y})) - h^0(\nn) \, , \label{eq-res1} \\
& &  \sum_{j=2}^{4} (-1)^j \cdot h^j(\Omega^3_{\tilde X}({\log \tilde Y}))
\stackrel{\scriptstyle{\eqref{eq-prel4}}}{=}
\sum_{j=2}^{4} ((-1)^j \cdot h^j(\Omega^3_{X}(Y))) + h^0(\nn) - \chi(\nn) \, . \label{eq-res2}
\end{eqnarray}

Finally,  from the equalities
$h^{3,2}(\tilde X) \stackrel{\scriptstyle{\eqref{eq-prel5}}}{=} h^{3,2}(X) =
h^{1,2}(X) = 0$,
the cohomology sequence associated to \eqref{poincare} breaks into shorter exact sequences and we obtain
\begin{eqnarray} \label{prawieh21}
h^{2,1}({\tilde Y})&=& h^{3,0}({\tilde X}) - h^{3,1}({\tilde X})
-
(h^0(\Omega^3_{\tilde X}({\log \tilde Y})) - h^1(\Omega^3_{\tilde X}({\log \tilde Y}))) +  h^{2,0}({\tilde Y}) = \nonumber \\
 &\stackrel{\scriptstyle{\eqref{eq-res1}}}{=}&  h^{3,0}(X) - h^{3,1}(X) -
h^0(\Omega^3_{X}({Y})) +  h^0(\nn ) + h^{2,0}({\tilde Y}) =
 h^{4,1}(X)  -
  h^{3,1}(X)  +  \nonumber \\
& &  -  h^4(\Omega^1_{X}({-Y})) +   h^0(K_X(2Y))
 -  h^0(K_X(Y))  - \mu_Y + \delta_Y  + h^{2,0}({\tilde Y}) \, .
\nonumber
\end{eqnarray}

By  \eqref{eq-res2} and similar argument we have
\bgroup\arraycolsep=3pt
\begin{eqnarray}  \label{prawieh11}
h^{2,2}({\tilde Y})  &=&
 h^{3,3}({\tilde X}) - h^{3,4}({\tilde X}) +
\sum_{j=2}^{4} (-1)^j \, h^j(\Omega^3_{\tilde X}({\log \tilde Y})) +  h^{2,3}({\tilde Y}) = \nonumber \\
&\stackrel{\scriptstyle{\eqref{eq-prel5}}}{=}&   h^{1,1}(X) + \mu_Y
 - h^{1,0}(X)  + \sum_{j=2}^{4} ((-1)^j \,  h^j(\Omega^3_{X}({Y}))) +  h^0(\nn) - \chi(\nn) +  h^{2,3}({\tilde Y})  \nonumber \\
  &\stackrel{\scriptstyle{\eqref{eq-chi}}}{=}&    h^{1,1}(X)
 + (\chi(\Omega^1_{X}({-Y})) - h^4(\Omega^1_{X}({-Y}))) - (\chi({\mathcal O}_{X}(-2Y)) - h^4({\mathcal O}_{X}(-2Y)))  \nonumber \\
 &&   -  h^{1,0}(X)  - h^1({\mathcal O}_{X}(-Y)) +   \mu_Y +  \delta_Y +   h^{2,3}({\tilde Y}) \, . \nonumber
\end{eqnarray}
\egroup

To complete the proof observe that
\eqref{eq-prel1} and the cohomology sequence associated to
\begin{equation} \label{eq-prel14}
0 \lra {\mathcal O}_{\tilde X}(-{\tilde Y}) \lra {\mathcal O}_{\tilde X}  \lra   {\mathcal O}_{\tilde Y}  \lra 0
\end{equation}
yield  (see also \cite[Ex.~II.8.8,~p.~190]{hart})
$$
h^{2,3}(\tilde{Y}) = h^{0,1}(\tilde{Y})
\stackrel{\scriptstyle{\eqref{eq-prel1}}}{=} h^{0,1}(X) - h^1({\mathcal
O}_{X}(-Y))
\, \, \mbox{ and } \, \,
 h^{2,0}(\tilde{Y}) =  h^{0,2}(X) \,.
$$
\end{proof}

We end this section with the study  of a classical example (see \cite{Clemens});
$k$-fold cyclic covers $Y$ of $\PP_3$ branched along a surface of degree $(d \cdot k)$.
 Recall that the assumptions of   neither \cite[Thm~1]{Cynk2}
nor  \cite[Thm~4.1]{rams-hab}
are fulfilled when we treat $Y$ as a hypersurface in the bundle
$\PP({\mathcal O}_{\PP_3} \oplus {\mathcal O}_{\PP_3}(d))$.
Below we check that  Thm~\ref{ade} works in that case.

\newcommand\tp{\tilde{\mathbb P}}

\begin{exmp}  \label{nagladkimidzie}
(c.f. \cite[Example~3.1]{rams-hab}) We fix integers $d \geq 2$ and $k \geq 2$. We consider
the $k$-fold cover $Y$  of $\PP_3$ branched along a surface of degree $d \cdot k$. We assume
that $Y$ has only  A-D-E singularities.
Let  $\tp = \PP({\mathcal E})$ with ${\mathcal E} := {\mathcal O}_{\PP_3} \oplus {\mathcal O}_{\PP_3}(d)$.
It is well-known that $Y$ can be considered as a hypersurface in $\tp$ and ${\mathcal O}_{\tp}(Y) = {\mathcal O}_{\tp}(k)$.
We claim that the pair $Y \subset \tp$ satisfies the assumptions
of Thm~\ref{ade}.

We maintain the notation of  \cite[Ex.~III.8.4]{hart}). Then we have $K_{\tp} =(\pi^{*}{\mathcal O}_{\PP_3}(d-4))(-2)$,
where $\pi : \tp \ra \PP_3$ stands for  the bundle projection. \\
At first we study the cohomologies of
$\Omega_{\tilde{\PP}/{\PP_3}}(-k)$.
We consider the exact sequence  \cite[Ex.~III.8.4.b]{hart} tensored with  ${\mathcal O}_{\tilde{\PP}}(-k)$:
\begin{equation} \label{eq-mik}
0 \lra \Omega_{\tilde{\PP}/{\PP_3}}(-k) \lra (\pi^{*}{\mathcal E})(-k-1)  \lra {\mathcal O}_{\tilde \PP}(-k)
\lra 0 \,  .
\end{equation}
In order to compute $ h^j({\mathcal O}_{\tilde{\PP}}(-k))$, observe that,
by \cite[Ex.~III.8.4.c]{hart}, we have
\begin{eqnarray*}
\pi_{*}((\pi^{*}{\mathcal O}_{\PP_3}(d-4))(k-2)) & = &
{\mathcal O}_{\PP_3}(d-4) \otimes  \mbox{S}^{k-2}{\mathcal E} = \bigoplus_{l=0}^{k-2} {\mathcal O}_{\PP_3}(ld + d - 4),  \\
\operatorname{R}^{j}\pi_{*}((\pi^{*}{\mathcal O}_{\PP_3}(d-4))(k-2)) &  = &  0
\quad \quad   \mbox{ for } j \geq 1, \,
\end{eqnarray*}
so the Leray spectral sequence and the Serre duality imply
\begin{equation} \label{eq-ch0}
h^4({\mathcal O}_{\tilde{\PP}}(-k)) = \sum_{l=0}^{k-2} \binom{ld+d-1}3
\mbox{ and }  h^j({\mathcal O}_{\tilde{\PP}}(-k)) = 0 \mbox{ for } j \leq 3 \, .
\end{equation}

We have $(\pi^{*}{\mathcal E})(-k-1) = {\mathcal O}_{\tp}(-k-1) \oplus  \pi^{*}{\mathcal O}_{\PP_3}(d)(-k-1)$.
As in \cite[(5.10)]{rams-hab} we use the Serre duality, \cite[Ex.~III.8.4.a]{hart}, the projection formula
and  the  Leray spectral sequence to show that
\begin{equation} \label{eq-ch3}
h^{1}(\pi^{*}{\mathcal E})(-k-1) = 1  \mbox{ and }  h^{j}(\pi^{*}{\mathcal E})(-k-1) = 0 \quad  \mbox{ for } j = 0, 2, 3.
\end{equation}

Now \eqref{eq-mik} and \eqref{eq-ch0}, \eqref{eq-ch3} give the
equalities
\begin{equation}  \label{koh1}
h^{1}(\Omega_{\tilde{\PP}/{\PP_3}}(-k)) = 1 \mbox{ and }  h^{j}( \Omega_{\tilde{\PP}/{\PP_3}}(-k)) = 0 \quad \mbox{ for } j = 0, 2, 3.
\end{equation}

In order to compute $h^{j}(\pi^{*} \Omega^1_{\PP_3}(-k))$, we consider  the pull-back of the Eu\-ler sequence under the map $\pi$
and tensor it with ${\mathcal O}_{\tilde{\PP}}(-k)$:
\begin{equation} \label{eq-ch4}
0 \lra (\pi^{*} \Omega^1_{\PP_3})(-k) \lra (\pi^{*}
{\mathcal O}_{\PP_3}(-1)^{\oplus4})(-k)
\lra {\mathcal O}_{\tilde{\PP}}(-k) \lra 0 \, .
\end{equation}
We use the Serre duality, \cite[Ex.~III.8.4.a]{hart} and the Leray spectral
sequence to see that
\begin{eqnarray*}
h^{4-j}(\pi^{*}{\mathcal O}_{\PP_{3}}(-1)(-k))  =
h^{j}(\pi^{*}{\mathcal O}_{\PP_{3}}(d-3)(k-2))
 =  h^{j}({\mathcal O}_{\PP_{3}}(d-3) \otimes \mbox{S}^{k-2} {\mathcal E}) \, .
\end{eqnarray*}
In this way we show that $ h^j( \pi^{*}({\mathcal O}_{{\PP_3}}(-1))(-k))= 0 \mbox{ for } j \leq3$.
The latter, combined with \eqref{eq-ch0} and \eqref{eq-ch4}, yields
\begin{equation} \label{eq-ch6}
  h^j(\pi^{*}\Omega^1_{\PP_3}(-k)  ) = 0  \mbox{ for } j \leq 3  \,  .
\end{equation}
Finally, we tensor  the exact sequence
\[
0 \lra \pi^{*} \Omega^1_{\PP_3}  \lra \Omega^1_{\tilde{\PP}} \lra  \Omega_{\tilde{\PP}/{\PP_3}}
\lra 0
\]
with ${\mathcal O}_{\tilde{\PP}}(-k)$, and apply \eqref{koh1}, \eqref{eq-ch6}  to see that
$h^j(\Omega^1_{\tp}(-Y))$ vanish for $j=2,3$, whereas $h^1(\Omega^1_{\tp}(-Y))=1$.

We use similar argument to show that  $h^2(\Omega^1_{\tp}) = 0$.
\end{exmp}

\section{Main result} \label{sect-singular}
\label{sec:singular}

In \cite{rams-hab} the ambient variety $X$ is assumed to be a projective normal Cohen-Macaulay fourfold.
In this section we study the question to what extent the formulae of  Thm~\ref{ade}
remain valid when  we allow the fourfold $X$ to be singular.

We assume $X$ to  be a {\sl four-dimensional normal complex variety}, so
the canonical (Weil) divisor  $K_X$ is well-defined (up to the linear equivalence).
Recall that  the map
$D \rightarrow  {\mathcal O}_X(D)$
gives  one-to-one correspondence between the linear equivalence classes of
Weil divisors and isomorphism classes of rank-$1$ reflexive sheaves on $X$ (see \cite[p.~281]{reid1}).
We put
$$
\OOmega :=  {\mathfrak j}_{*}\Omega^3_{\mbox{\tiny reg}(X)},
$$
 where ${\mathfrak j}:\mbox{reg}(X) \rightarrow X$ stands for the inclusion.

Let $Y \subset X$ be a  a hypersurface with A-D-E singularities such that
\begin{equation} \label{eq-omijaosobliwosci}
\mbox{sing}(X) \cap Y = \emptyset \, .
\end{equation}
We maintain the notation of the previous section. In particular,
$\tilde{Y}$ is given by \eqref{eq-prel6}.
Observe that, by the assumption \eqref{eq-omijaosobliwosci}, we can consider the pullbacks $\ppi^{\ast}K_X$, $\ppi^{\ast}Y$.
 Obviously, there exist unique positive integers $k_l, m_l$ satisfying the conditions of Def.~\ref{def-ideal}.
 In particular, we have  the linear equivalence \eqref{eq-cof-id}.

We define {\sl the ideal} $\II_{Y}$ (resp. {\sl the defect of
}$Y$) by the equality \eqref{eq-prel1ipol} (resp.
\eqref{eq-def}). \\
Let $\Oomega^3_{\tilde X}(\ast \, \tilde{Y}) := \lim\limits_{\substack{\longrightarrow\\k}}
 \, \Oomega^3_{\tilde X}(k \tilde{Y})$, and let  $V$ be an open subset of $\tilde{X}$. We put
\begin{eqnarray*}
\Gamma(V,\Oomega^3_{\tilde X}(\log \tilde{Y})) &:=& \{ \alpha \in  \Gamma(V,\Oomega^3_{\tilde X}(\ast \, \tilde{Y}));
\alpha|_{\mbox{\tiny reg}(X) \cap V} \mbox{ and } \mbox{d}(\alpha|_{\mbox{\tiny reg}(X) \cap V})
 \mbox{ have at most  \hspace*{4ex}} \\
& &  \mbox{  \hspace*{40ex} simple  poles along } \tilde{Y} \},
\end{eqnarray*}
The assumption \eqref{eq-omijaosobliwosci}, combined with
\eqref{poincare}, \eqref{resolution}, implies that the sequences
\begin{equation} \label{eq-exsbis}
0 \lra  \Oomega^3_{\tilde X} \lra  \Oomega^3_{\tilde X}(\log \tilde{Y}) \lra  \Omega^2_{\tilde Y} \lra  0   \, ,
\end{equation}
\begin{equation} \label{eq-exstertio}
0 \lra  \Oomega^3_{\tilde X}(\log \tilde{Y}) \lra  \Oomega^3_{\tilde X}(\tilde{Y}) \lra
{\mathcal O}_{\tilde{X}}(2 \tilde{Y} + K_{\tilde{X}}) / {\mathcal O}_{\tilde{X}}(\tilde{Y} + K_{\tilde{X}}) \lra 0 \,
\end{equation}
are exact.

Now we are in position to  prove a more general version of  \eqref{eq-prell}, \eqref{eq-prel4}, \eqref{eq-prel5}.
\begin{lem} \label{lem-osobliwe}
  We have the following equalities: \\
a) $h^{3}(\Oomega^3_{\tilde X}) = h^{3}(\OOmega) + \mu_Y$ and
$h^{i}(\Oomega^3_{\tilde X}) = h^{i}(\OOmega) \mbox{ for } i \not=3,$ \\
b) $h^{j}(\Oomega^3_{{\tilde X}}(\tilde{Y}))  = h^{j}(\Oomega_X^3(Y)) $  for all $j$, \\
c) $h^{j}({\mathcal O}_{\tilde{X}}( \tilde{Y} + K_{\tilde{X}})  ) =  h^{j}({\mathcal O}_X(Y + K_X))$ for all $j$, \\
d) $ \chi({\mathcal O}_{\tilde{X}}(2 \tilde{Y} + K_{\tilde{X}})) = \chi({\mathcal O}_{X}(2 Y + K_X)) - \mu_Y$, \\
e)  $h^{j}({\mathcal O}_{\tilde{X}}( - \tilde{Y})) =  h^{j}({\mathcal O}_X(-Y))$ for all $j$.
\end{lem}
\begin{proof}
Let $E$ stand for the exceptional divisor of the blow-up $\sigma_1$.
By direct computation
\begin{equation} \label{eq-pullback}
(\sigma_1)^{\ast}(\Oomega_X^3) = \Oomega_{\tilde X^1}^3(\log E)(-3E) \, .
\end{equation}
Obviously, we have $(\sigma_1)_{\ast}{\mathcal O}_{\tilde{X}^1}(kE) = {\mathcal O}_{X}$ for  $k = 1, 2,3$.
We follow the proof of \cite[(2.5)]{rams-hab} to show that
\begin{equation} \label{eq-ideal}
R^{j}(\sigma_1)_{\ast} {\mathcal O}_{\tilde{X}^1}(lE) = 0 \mbox{
where } j > 0, \,  l \leq 3.
\end{equation}
Therefore, from the projection
formula (see \cite[(2.8)]{rams-hab}), we get for $k=1,2,3$
\begin{equation} \label{eq-dirh}
(\sigma_1)_{\ast}(\sigma_1^{\ast}(\Oomega_X^3)(kE)) = \Oomega_X^3 \mbox{ and }
R^{j}(\sigma_1)_{\ast}( \sigma_1^{\ast}(\Oomega_X^3)(lE)) = 0 \mbox{ where }j > 0 \mbox{, } l \leq 3.
\end{equation}

\noindent
{\sl a)} The projection formula, combined with \eqref{eq-pullback} and \eqref{eq-dirh}, yields
\begin{equation} \label{eq-dir-log}
(\sigma_1)_{\ast}(\Oomega_{\tilde X^1}^3(\log E)) = \Oomega_{X}^3, \, \,  \, \, \,
 R^{j}(\sigma_1)_{\ast}(\Oomega_{\tilde X^1}^3(\log E)) =  0 \mbox{ for } j > 0 \, .
\end{equation}
We consider cohomology sequence associated to the exact sequence (see \cite[p.~444]{pet-st})
\begin{equation} \label{eq-logE}
0 \lra  \Oomega^3_{\tilde X^1} \lra  \Oomega^3_{\tilde X^1}(\log E) \lra  \Omega^2_{E} \lra 0
\end{equation}
to obtain the equalities
\begin{equation} \label{eq-h4}
h^4(\Oomega^3_{\tilde X^1}) = h^4(\Oomega^3_{\tilde X^1}(\log E)) = h^{4}(\Oomega_{X}^3),
\end{equation}
where the latter results from the Leray spectral sequence and \eqref{eq-dir-log}.

Observe that \eqref{eq-pullback}, combined with \cite[2.3~Property~c]{esnview}
for $D = D_1 = E$,
yields the exact sequence:
\begin{equation}  \label{eq-esv}
0 \lra \sigma_1^{\ast}(\Oomega_X^3)(2E) \lra  \Oomega^3_{{\tilde X}^1} \lra  \Omega^3_{E} \lra  0   \, .
\end{equation}
We consider the direct image of \eqref{eq-esv} under $\sigma_1$.
The centers of $\sigma_1$ are   smooth points on $X$, so
$R^{3}(\sigma_1)_{\ast} \Omega^3_{E} =  \CC^{\nu_{0}}$  is the sky-scraper
sheaf with stalks $\CC$ in the centers of the blow-up and
$R^{i}(\sigma_1)_{\ast} \Omega^3_{E}$  vanish for $i = 0,1,2$.
The latter, combined with \eqref{eq-dirh}, implies
\begin{equation*}
(\sigma_1)_{\ast}\Oomega_{\tilde X^1}^3 = \Oomega_{X}^3, \, \,
 R^{i}(\sigma_1)_{\ast} \Oomega_{\tilde X^1}^3 =  0 \mbox{ for } i = 1,2 \, ,
\mbox{ and }
 R^{3}(\sigma_1)_{\ast} \Oomega_{\tilde X^1}^3 = \CC^{\nu_{0}}.
\end{equation*}
From the Leray
  spectral sequence (see e.g. \cite[Example 1.D]{mccleary}), we obtain the equalities
$h^{i}(\Oomega_{X}^3) = h^{i}(\Oomega_{\tilde X^1}^3)$  for $i = 0, 1, 2$ and the exact sequence
\begin{equation*} \label{eq-dlugi}
0
 \ra
H^{3}(\Oomega_{X}^3)
\lra
H^{3}(\Oomega_{\tilde X^1}^3)
\lra
\CC^{\nu_{0}}
\lra
H^{4}( \Oomega_{X}^3)
\lra
H^{4}(\Oomega_{\tilde X^1}^3)
\lra
0 \, .
\end{equation*}
Therefore, \eqref{eq-h4} yields that
$h^3(\Oomega^3_{\tilde X^1}) =  h^{3}(\Oomega_{X}^3) + \nu_0.$

To complete the proof of the part a) of the lemma proceed by induction on the number of blow-ups in \eqref{eq-prel6}.

\noindent
{\sl b)} Observe that, as in the proof of a), it suffices to show the equalities
\begin{equation} \label{eq-bla}
h^{i}(\Oomega_{\tilde{X}^1}^3(\tilde{Y}^1))  =  h^{i}(\Oomega_X^3(Y))  \mbox{ for } i \geq 0 \, .
\end{equation}
 We claim that
\begin{equation} \label{eq-dirh2}
(\sigma_1)_{\ast}(\Oomega_{{\tilde X}^1}^3(-2E)) = \Oomega_X^3 \quad\mbox{ and
}\quad
R^{j}(\sigma_1)_{\ast}(\Oomega_{{\tilde X}^1}^3(-2E)) = 0 \mbox{  for }j > 0.
\end{equation}
Indeed, we tensor the exact sequence \eqref{eq-logE} with the (locally free) sheaf ${\mathcal O}_{\tilde{X}^1}(-2E)$
to obtain:
\begin{equation} \label{eq-ex2}
0 \lra  \Oomega^3_{{\tilde X}^1}(-2E) \lra   \sigma_1^{\ast}(\Oomega_{{\tilde X}^1}^3)(E) \lra   \Omega^2_{E}(-2E) \lra  0   \, .
\end{equation}
We consider the direct image of \eqref{eq-ex2} under $\sigma_1$.
Since
$\Omega^2_{E}(-2E) = \Omega^2_{\PP_3}(2)$, we have the vanishings
$R^{l}(\sigma_1)_{\ast}\Omega^2_{E}(-2E) = 0$ for $l \geq 0$. Now \eqref{eq-dirh2}  results immediately from
\eqref{eq-dirh}.   \\
The equality
${\mathcal O}_{\tilde{X}^1}(\tilde{Y}^1)  =
(\sigma_1)^{\ast}{\mathcal O}_{X}(Y) \otimes {\mathcal
O}_{\tilde{X}^1}(-2E)$
(see also \cite[(2.2)]{rams-hab}),
combined with the projection formula  and \eqref{eq-dirh2},
yields
\begin{equation*}
(\sigma_1)_{\ast}(\Oomega_{\tilde{X}^1}^3(\tilde{Y}^1)) =
\Oomega_X^3(Y) \mbox{ and }
R^j(\sigma_1)_{\ast}(\Oomega_{\tilde{X}^1}^3(\tilde{Y}^1)) = 0 \mbox{ for } j > 0.
\end{equation*}
Thus the Leray spectral sequence implies \eqref{eq-bla}.

\noindent
{\sl c)} Let $k = 1,2$. From the assumption \eqref{eq-omijaosobliwosci}, we have
\begin{equation} \label{eq-hdirvan}
{\mathcal O}_{\tilde{X}^1}(k \tilde{Y}^1 + K_{\tilde{X}^1}) = (\sigma_1)^{\ast}{\mathcal O}_{X}(k Y  +  K_X)
 \otimes  {\mathcal O}_{\tilde{X}^1}((3-2k) E).
\end{equation}
We apply the projection formula
 and \eqref{eq-ideal} to show that
\begin{equation} \label{eq-nareszcie}
(\sigma_1)_{\ast}{\mathcal O}_{\tilde{X}^{1}}(\tilde{Y}^{1} +
K_{\tilde{X}^{1}}) = {\mathcal O}_{X}(Y  +  K_X)
\, \mbox{ and } \,
R^p (\sigma_1)_{\ast}{\mathcal O}_{\tilde{X}^1}(k \tilde{Y}^1 + K_{\tilde{X}^1}) = 0  \mbox{ for } p > 0.
\end{equation}
The Leray spectral sequence and \eqref{eq-nareszcie} yield the
equalities
$$
h^{j}({\mathcal O}_{\tilde{X}^1}( \tilde{Y}^1 + K_{\tilde{X}^1})  ) =  h^{j}({\mathcal O}_{X}(Y + K_X)) \mbox{ for all } j.
$$

\noindent
{\sl d)} Let $i = 1, \ldots, n$. Reasoning as in
\eqref{eq-hdirvan}, \eqref{eq-nareszcie}
we obtain
\begin{eqnarray*}
{(\sigma_{i})}_{*}O_{\tilde{X}^{i}}(K_{\tilde{X}^{i}}+2\tilde{Y}^
{i
} )&=&O_{\tilde{X}^{i-1}}(K_{\tilde{X}^{i-1}}+2\tilde{Y}^{i-1}
)\otimes \mathcal J_{i}  \, ,  \\
{R^{p}(\sigma_{i})}_{*}O_{\tilde{X}^{i}}(K_{\tilde{X}^{i}}
+2\tilde {
Y }^{ i} )&=&0\quad\text{ for }p>0,
\end{eqnarray*}
where $\mathcal J_{i}$ is the (reduced) ideal of the center of
the blow--up $\sigma_{i}$. From the Leray spectral sequence we
get
\[\chi(O_{\tilde{X}^{i}}(K_{\tilde{X}^{i}}+2\tilde{Y}^{i
}))=\chi(O_{\tilde{X}^{i-1}}(K_{\tilde{X}^{i-1}}+2\tilde{Y}^{i-1
})\otimes\mathcal  J_{i}).\]
We have the exact sequence
\[
0\ra    O_{\tilde{X}^{i-1}}(K_{\tilde{X}^{i-1}}+2\tilde{Y}^{i-1
}) \otimes \mathcal J_{i} \ra
O_{\tilde{X}^{i-1}}(K_{\tilde{X}^{i-1}}+2\tilde{Y}^{i-1
})\ra \mathcal S\ra0,
\]
where $\mathcal S = \CC^{\nu_{i-1}}$ is the sky--scraper sheaf
with stalks $\CC$   over the points in
$\sing(\tilde{Y}^{i-1})$. Hence
\[\chi(O_{\tilde{X}^{i-1}}(K_{\tilde{X}^{i-1}}+2\tilde{Y}^{i-1
})\otimes\mathcal
J_{i})=\chi(O_{\tilde{X}^{i-1}}(K_{\tilde{X}^{i-1}}+2\tilde{Y}^{
i-1}))-\nu_{i}.\]

\noindent
{\sl e)} The equality
${\mathcal O}_{\tilde{X}^1}(- \tilde{Y}^1)  = (\sigma_1)^{\ast}{\mathcal O}_{X}(-Y) \otimes {\mathcal O}_{\tilde{X}^1}(2E)$,
combined with \eqref{eq-ideal}, the projection formula and the Leray spectral sequence, implies that
$$
h^{j}({\mathcal O}_{\tilde{X}^1}(- \tilde{Y}^1)) =  h^{j}({\mathcal O}_{X}(- Y)) \mbox{ for all } j.
$$
\end{proof}

Now we are in position to prove
\begin{thm} \label{ade-osobliwe}
Let $X$ be a four-dimensional normal complete variety,
and let $Y \subset X$ be a hypersurface with
A-D-E singularities such that $\mbox{sing}(X) \cap Y =
\emptyset$.
If
$$
 h^1({\mathcal O}_{X}(Y + K_{X})) = h^2(\Oomega^3_{X}) = h^1(\Oomega^3_{X}(Y) =
0 \quad\mbox{ and }\quad
h^{2}({\mathcal O}_{X}(-Y))= h^{3}({\mathcal O}_{X}(-Y)) = 0 \, ,
$$
then the following equalities hold
\bgroup\arraycolsep=1pt
\begin{eqnarray*}
 h^{1,1}({\tilde Y})&=& h^{3}(\OOmega) +
 (\chi(\OOmega(Y)) - h^0(\OOmega(Y))) + (\chi({\mathcal O}_{X}(Y + K_{X}))   - h^0({\mathcal O}_{X}(Y + K_{X}))) + \nonumber \\
 & &  + h^1({\mathcal O}_X) - h^4(\OOmega) - (\chi({\mathcal O}_{X}(2Y + K_{X})) - h^0({\mathcal O}_{X}(2Y + K_{X})))
-  h^1({\mathcal O}_{X}(-Y)) \nonumber \\
 & &  +  \mu_Y +  \delta_Y  \, , \nonumber \\
h^{1,2}({\tilde Y})&=&  h^0(\OOmega) + h^2({\mathcal O}_X) + h^0({\mathcal O}_{X}(2Y + K_{X})) -
  h^1(\OOmega)  - h^0({\mathcal O}_{X}(Y + K_{X}))  + \nonumber \\
 & &  -  h^0(\OOmega(Y))  - \mu_Y + \delta_Y \, , \nonumber
\end{eqnarray*}
\egroup
where $\delta_Y$ (resp. $\mu_Y$) is the defect
(resp. the number of singularities and infinitely near singularities) of $Y$.
\end{thm}
\begin{proof}
Lemma~\ref{lem-osobliwe} and the exact sequences  \eqref{eq-exsbis}, \eqref{eq-exstertio}
enable us to repeat the proof of Thm~\ref{ade}: \\
We define the quotient sheaf $\nn$ (see \eqref{eq-exstertio}) and use the vanishing
of $h^1({\mathcal O}_{X}(Y + K_{X}))$, combined with Lemma~\ref{lem-osobliwe}.d,  to compute $h^0(\nn)$. \\
Then we apply Lemma~\ref{lem-osobliwe} to study the cohomology sequences given by  \eqref{eq-exsbis} and \eqref{eq-exstertio}. \\
Finally, we use   Lemma~\ref{lem-osobliwe}.e and the sequence \eqref{eq-prel14} to calculate the numbers $h^{2,3}(\tilde{Y})$, $h^{2,0}(\tilde{Y})$ (observe that we blow up smooth points on $X$, so the Leray spectral sequence immediately yields
$h^j({\mathcal O}_X) = h^j({\mathcal O}_{\tilde X})$  for each $j$).
\end{proof}

\noindent
Remarks:
1. In certain applications $h^{1,0}(\tilde{Y})$, $h^{2,0}(\tilde{Y})$ are known. In such case, only the first three vanishings
that we assume in Thm~\ref{ade-osobliwe}
are needed to find all the Hodge numbers of $\tilde{Y}$. For example, if we assume $X$ to be a (normal) toric variety, then
$h^i({\mathcal O}_X) = 0$ for $i > 0$ (see  \cite[Corollary~2.8]{oda}), so we have
$h^i({\mathcal O}_{\tilde Y}) = h^{i+1}({\mathcal O}_{X}(-Y))$. \\
2. If the Serre duality holds on $X$ with ${\mathcal O}_X(K_X)$ as the dualizing sheaf , e.g. $X$ is
a projective normal
Cohen-Macaulay fourfold, the assumption on $h^{3}({\mathcal O}_{X}(-Y))$ is abundant. \\
3.  The Zariski sheaf of germs of $3$-forms  $\OOmega$ satisfies the condition
$\OOmega |_{\mbox{\tiny reg}(X)} =  \Omega^3_{\mbox{\tiny reg}(X)}$.
Obviously, our considerations remain true if we replace $\OOmega$ with any sheaf with the above property.

\vspace*{3mm}

Recall that
one can use \cite[Thm~4.1]{rams-hab} to compute the Hodge numbers
of ample hypersurfaces in
projective toric {\sl simplicial} varieties (see
\cite[Cor.~4.2]{rams-hab}).
By Danilov's spectral sequence
Thm~\ref{ade-osobliwe} works for all projective (normal) toric varieties $X$ as the following corollary shows:

\begin{cor} \label{toric}
Let $X$ be a complete toric fourfold,  and let $Y \subset X$ be a hypersurface with A-D-E singularities
 such that $\sing(X) \cap Y = \emptyset$. If $\mo_{X}(Y)$ is ample, then
\begin{eqnarray*}
 h^{1,1}({\tilde Y})&=& h^{3}(\OOmega)  +  \mu_Y +  \delta_Y  \, , \nonumber \\
h^{1,2}({\tilde Y})&=&  h^0({\mathcal O}_{X}(2Y + K_{X}))   - h^0({\mathcal O}_{X}(Y + K_{X})) -  h^0(\OOmega(Y))  - \mu_Y + \delta_Y \, . \nonumber
\end{eqnarray*}
\end{cor}
\begin{proof}
 By Danilov's Spectral sequence (see \cite[p.~133]{oda}) we have  $h^j(\Oomega^3_{X}) = 0$ for $j < 3$
and $h^4(\Oomega^3_{X}) = 0$.
\cite[Bott's Vanishing,~p.~130]{oda} (see also
\cite[Cor.~1.3]{fujino})
implies that $h^i(\OOmega(Y))$ and $h^i({\mathcal O}_{X}(kY +
K_{X}))$ vanish for $k=1,2$ and $i > 0$.
We use the latter vanishing and the Serre duality to show that $h^i({\mathcal
O}_{X}(-Y)) = 0$ for  $i > 0$.
\end{proof}

\noindent
A nef Cartier divisor $D$ on a complete simplicial toric variety is called $l$-semiample iff
its Kodaira-Itaka dimension equals $l$. For such divisors one can use \cite[Cor.~2.7]{Mavlyutov} to compute
$h^i(\Oomega^3_{X}(Y))$ and check whether the assumptions of  Thm~\ref{ade-osobliwe} are fulfilled.
In particular, we have

\begin{exmp} \label{1-ample}
Let $X$ be a complete simplicial toric variety, and let $Y \subset X$ be a hypersurface with A-D-E singularities
satisfying \eqref{eq-omijaosobliwosci} and such that ${\mathcal O}_{X}(Y)$ is $1$-semiample. Then, by \cite[Thm~2.4]{Mavlyutov},
$Y$ satisfies the assumptions of Thm~\ref{ade-osobliwe}.
\end{exmp}

\noindent
Remark (c.f.  Ex.~\ref{nagladkimidzie}): One can easily see that the
$k$-fold solid $Y \subset \PP({\mathcal O}_{\PP_3} \oplus {\mathcal O}_{\PP_3}(d))$
defines a nef and $4$-ample divisor. It is not ample because  $h^1(\Omega^1_{\tp}(-Y)) \neq 0$ (see \cite[Cor.~6.4]{esnview}).

\VVsect

\vspace*{3ex}
\noindent
S{\l}awomir Cynk \\
Institute of Mathematics,
Jagiellonian University,
ul. {\L}ojasiewicza 6,
30-348 Krak\'ow,
Poland \\[1mm]
Institute of Mathematics of the
Polish Academy of Sciences, ul. \'Sniadeckich 8, 00-956 Warszawa, Poland \\
e-mail: Slawomir.Cynk@im.uj.edu.pl  \\[2mm]
S{\l}awomir Rams \\
 Department Mathematik,
 Universit\"at Erlangen-N\"urnberg,
 Bismarckstra\ss e 1 1/2,
 D-91054 Erlangen,
 Germany \\
 and \\
Institute of Mathematics,
Jagiellonian University,
ul.  {\L}ojasiewicza 6,
30-348 Krak\'ow,
Poland  \\
e-mail: Slawomir.Rams@im.uj.edu.pl

\parindent=0cm

\begin{thebibliography}{99}
\bibitem{bat93}   V.~Batyrev, D.~A.~Cox,
\emph{On the Hodge structure of projective hypersurfaces in toric varieties.}
Duke Math. J. \textbf{75} (1994), 293--338.
\bibitem{Chelt}  I.~Cheltsov, \emph{On factoriality of nodal
threefolds.}
 J. Alg. Geom. \textbf{14} (2005), 663--690.
\bibitem{Clemens} C.~H.~Clemens, \emph{Double solids.} Adv. in Math. \textbf{47} (1983), 107--230.
\bibitem{Cynk2} S.~Cynk, \emph{Defect of a nodal hypersurface.}
Manuscripta Math. \textbf{104} (2001), 325--331.
\bibitem{Cynk3} S. Cynk,  \emph{Cohomologies of a double covering of a non-singular algebraic $3$-fold.}
Math. Z. \textbf{240} (2002), 731--743.
\bibitem{Cynk4} S. Cynk, D. van Straten,  \emph{Infinitesimal deformations of double covers of smooth
algebraic varieties.}  Math. Nach. \textbf{279} (2006), 716--726.
\bibitem{digennaro} V.~Di Gennaro, D.~Franco, \emph{ Factoriality
and N\'eron-Severi groups.}
Commun. Contemp. Math.  \textbf{10} (2008),  745--764.
\bibitem{Dimca} A.~Dimca, \emph{Betti numbers of hyperplanes and defects of
linear systems}, Duke Math. Jour. \textbf{60} (1990),285--294.
\bibitem{Dimca2} A.~Dimca, \emph{Singularities and topology of hypersurfaces.} Springer 1992.
\bibitem{esnview}  H.~Esnault, E.~Viehweg, \emph{Lectures on vanishing theorems.} Birkh\"auser 1992.
\bibitem{fujino} O.~Fujino, \emph{Multiplication maps and vanishing theorems for Toric varieties.}
Math. Z. \textbf{257} (2007), 631-641.
\bibitem{hart} R.~Hartshorne, \emph{Algebraic geometry.} Springer 1977.
\bibitem{hulek}  K.~Hulek, R.~Kloosterman,  \emph{Calculating the Mordell-Weil rank of elliptic threefolds and the cohomology of singular hypersurfaces.} Preprint available at arXiv: math/0806.2025, 2008.
\bibitem{Kollar} J.~Koll\'ar, S.~Mori, \emph{Birational geometry of algebraic varieties.}  Cambridge University Press 1998.
\bibitem{Mavlyutov} A.~Mavlyutov, \emph{Cohomology  of rational forms and a vanishing theorem on toric varieties.}
J.~Reine Angew. Math. \textbf{61} (2008), 45--58.
\bibitem{mccleary}
J.~McCleary, \emph{User's guide to spectral sequences.} Cambridge University Press 2001.
\bibitem{oda}  T.~Oda,
\emph{Convex bodies and algebraic geometry. An introduction to the theory of toric varieties.}
 Springer 1988.
\bibitem{pet-st} C.~Peters, J.~Steenbrink, \emph{Infinitesimal variations of Hodge structure and the generic Torelli problem
for projective hypersurfaces.}  Classification of
algebraic and analytic manifolds (Katata, 1982), 399--463, Progr. Math. 39, Birkh\"auser 1983.
\bibitem{pet-st-book} C.~Peters, J.~Steenbrink, \emph{Mixed Hodge Structures.} Springer 2008
\bibitem{cheltsov2} V.~V.~Przhiyalkovskii, I.~Cheltsov, K.~A.~~Shramov,
 \emph{ Hyperelliptic and trigonal Fano threefolds.} (Russian)
Izv. Ross. Akad. Nauk Ser. Mat.  \textbf{69}
(2005), 145--204.
\bibitem{rams-hab} S.~Rams,
\emph{Defect and Hodge numbers of hypersurfaces.}
 Adv. Geom. \textbf{8} (2008), 257-288.
\bibitem{reid1} M.~Reid, \emph{Canonical 3-folds.}
Journees de geometrie algebrique, Angers/France 1979, 273-310.
\bibitem{vanstrquint}  D. van Straten,  \emph{A quintic hypersurface in $\PP_4$ with 130 nodes},
Topology~\textbf{32} (1993), 857-864.
\bibitem{vS} D. van Straten,
 \emph{Gutachten \"uber die Habilitationsschrift ``Defect and Hodge Numbers of Hypersurfaces''.} Mainz 2006.

\end{thebibliography}
\end{document}